\documentclass{amsart}
\usepackage{amssymb, amsthm, amscd, amsmath}

\newtheorem{theorem}{Theorem}[section]
\newtheorem{lemma}[theorem]{Lemma}
\newtheorem{prop}[theorem]{Proposition}

\theoremstyle{definition}

\theoremstyle{remark}
\newtheorem{remark}[theorem]{Remark}

\def\image{\mathop{\mathrm{Image}}}
\def\height{\mathop{\mathrm{height}}}

\begin{document}

\title{$p$-torsion elements in local cohomology modules}

\author{Anurag K. Singh}
              
\address{Department of Mathematics, University of Utah, 155 S. 1400
E., Salt Lake City, \quad UT 84112, USA \newline
E-mail: {\tt singh@math.utah.edu}}

\begin{abstract}  
For every prime integer $p$, M.~Hochster conjectured the existence of certain
$p$-torsion elements in a local cohomology module over a regular ring of mixed
characteristic. We show that Hochster's conjecture is false. We next construct 
an example where a local cohomology module over a hypersurface has $p$-torsion
elements for every prime integer $p$, and consequently has infinitely many
associated prime ideals.  
\end{abstract}

\maketitle

For a commutative Noetherian ring $R$ and an ideal ${\mathfrak a} \subset R$,
the finiteness properties of the local cohomology modules $H_{\mathfrak
a}^i(R)$ have been studied by various authors. In this paper we focus on the
following question raised by C.~Huneke \cite[Problem 4]{Hu}: if $M$ is a
finitely generated $R$-module, is the number of associated primes ideals of
$H_{\mathfrak a}^i(M)$ always finite?

In the case that the ring $R$ is regular and contains a field of prime
characteristic $p > 0$, Huneke and Sharp showed in \cite{HS} that the set of
associated prime ideals of $H_{\mathfrak a}^i(R)$ is finite. If $R$ is a
regular {\it local}\/ ring containing a field of characteristic zero,
G.~Lyubeznik showed that $H_{\mathfrak a}^i(R)$ has only finitely many
associated prime ideals, see \cite{Ly1} and also \cite{Ly2, Ly3}. Recently
Lyubeznik has also proved this result for unramified regular local rings of
mixed characteristic, \cite{Ly4}. Our computations here support Lyubeznik's
conjecture \cite[Remark 3.7 (iii)]{Ly1} that local cohomology modules over all
regular rings have only finitely many associated prime ideals. 

In Section \ref{elem} we construct an example of a hypersurface $R$ such that
the local cohomology module $H_{\mathfrak a}^3(R)$ has $p$-torsion elements for
every prime integer $p$, and consequently has infinitely many associated prime
ideals.  

For some of the other work related to this question, we refer the reader to
the papers \cite{BL, BRS, He}.

\section{Hochster's conjecture}

Consider the polynomial ring over the integers $R={\mathbb Z}[u,v,w,x,y,z]$
where ${\mathfrak a}$ is the ideal generated by the size two minors of the
matrix 
$$
M = \begin{pmatrix}  
u & v & w \\  
x & y & z \\  
\end{pmatrix}, 
$$
i.e., ${\mathfrak a} = (\Delta_1, \Delta_2, \Delta_3)R$ where $\Delta_1 =
vz-wy$, $\Delta_2 = wx-uz$, and $\Delta_3 = uy-vx$.  M.~Hochster conjectured
that for every prime integer $p$, there exist certain $p$-torsion elements in
the local cohomology module $H_{\mathfrak a}^3(R)$, and consequently that
$H_{\mathfrak a}^3(R)$ has infinitely many associated prime ideals. We first
describe the construction of these elements. For an arbitrary prime integer 
$p$, consider the short exact sequence
$$
\begin{CD}
0 @>>> R @>p>> R @>>> R/pR @>>> 0,
\end{CD}
$$
which induces the following exact sequence of local cohomology modules:
$$
0 \to H_{\mathfrak a}^2(R) \to H_{\mathfrak a}^2(R) 
\to H_{\mathfrak a}^2(R/pR) \to H_{\mathfrak a}^3(R) 
\to H_{\mathfrak a}^3(R) \to 0. \qquad (*)
$$
It can be shown that the module $H_{\mathfrak a}^3(R/pR)$ is zero using a 
result of Peskine and Szpiro, see Proposition \ref{PS} below. In the ring $R$ 
we have the equation
$$
u\Delta_1 + v\Delta_2 + w\Delta_3 = 0,
$$ 
arising from the determinant of the matrix
$$
\begin{pmatrix} 
u & v & w \\ 
u & v & w \\ 
x & y & z \\ 
\end{pmatrix}. 
$$
Let $q=p^e$ where $e$ is a positive integer. Taking the $q$\,th power of the 
above equation, we see that 
$$
(u\Delta_1)^q + (v\Delta_2)^q + (w\Delta_3)^q \equiv 0 \mod p, \qquad (**)
$$ 
and this yields a relation on the elements $\overline{\Delta}_1^q, \
\overline{\Delta}_2^q, \ \overline{\Delta}_3^q \in R/pR$, where 
$\overline{\Delta}_i$ denotes the image of $\Delta_i$ in $R/pR$. This relation
may be viewed as an element $\mu_q \in H_{\mathfrak a}^2(R/pR)$. Hochster 
conjectured that for every prime integer there exists a choice of 
$q=p^e$ such that
$$
\mu_q \notin \image \Big(H_{\mathfrak a}^2(R) \to H_{\mathfrak a}^2(R/pR)\Big),
$$
and consequently that the image of $\mu_q$ under the connecting homomorphism in
the exact sequence $(*)$ is a nonzero $p$-torsion element of 
$H_{\mathfrak a}^3(R)$. 

We note an equivalent form of Hochster's conjecture which is convenient to work
with. Recall that the module $H_{\mathfrak a}^3(R)$  may be viewed as the
direct limit
$$ 
\varinjlim \frac{R}{(\Delta_1^k, \Delta_2^k, \Delta_3^k)R}
$$ 
where the maps in the direct limit system are induced by multiplication by the 
element $\Delta_1 \Delta_2 \Delta_3$. The equation $(**)$ shows that
$$
\lambda_q = \frac{(u\Delta_1)^q + (v\Delta_2)^q + (w\Delta_3)^q}{p}
$$ 
has integer coefficients, i.e., that $\lambda_q$ is an element of the ring $R$. 
Let
$$
\eta_q = [\lambda_q + (\Delta_1^q, \ \Delta_2^q, \ \Delta_3^q)R] \ \in \ 
H_{\mathfrak a}^3(R). 
$$

\begin{lemma}
With the notation as above, the following statements are equivalent:
\item $(1)$ $\mu_q \notin \image \Big(H_{\mathfrak a}^2(R) 
\to H_{\mathfrak a}^2(R/pR)\Big)$,
\item $(2)$ $\eta_q$ is a nonzero $p$-torsion element of $H_{\mathfrak a}^3(R)$. 
\end{lemma}

\begin{proof}
$(2) \implies (1)$ If $\mu_q \in \image \Big(H_{\mathfrak a}^2(R) \to 
H_{\mathfrak a}^2(R/pR)\Big)$, the relation $(\overline{u}^q, \overline{v}^q, 
\overline{w}^q)$ on the elements $\overline{\Delta}_1^q, \ \overline{\Delta}_2^q, 
\ \overline{\Delta}_3^q\in R/pR$ \, lifts to a relation in $H_{\mathfrak
a}^2(R)$, i.e., there exists an integer $k$ and elements $\alpha_i \in R$ such 
that $\alpha_1 \Delta_1^{q+k}+\alpha_2 \Delta_2^{q+k}+\alpha_3 \Delta_3^{q+k}=0$ 
and
$$
(u^q\Delta_2^k\Delta_3^k, \ v^q\Delta_3^k\Delta_1^k, \ w^q\Delta_1^k\Delta_2^k)
\equiv (\alpha_1, \ \alpha_2, \ \alpha_3) \mod p.
$$
Hence we have
\begin{align*}
& (u^q\Delta_2^k\Delta_3^k -\alpha_1)\Delta_1^{q+k} 
+(v^q\Delta_3^k\Delta_1^k -\alpha_2)\Delta_2^{q+k}
+(w^q\Delta_1^k\Delta_2^k -\alpha_3)\Delta_3^{q+k} \\
&=\Big((u\Delta_1)^q+(v\Delta_2)^q+(w\Delta_3)^q\Big)(\Delta_1\Delta_2\Delta_3)^k
\in (p\Delta_1^{q+k}, \ p\Delta_2^{q+k}, \ p\Delta_3^{q+k} )R,
\end{align*}
and so $\lambda_q (\Delta_1 \Delta_2 \Delta_3)^{k} \in
(\Delta_1^{q+k}, \ \Delta_2^{q+k}, \ \Delta_3^{q+k})R$ and $\eta_q=0$.

The argument that $(1) \implies (2)$ is similar.
\end{proof}

By the above lemma, an equivalent form of Hochster's conjecture is that for
every prime integer $p$, there is an integer $q=p^e$ such that $\eta_q$ is a
nonzero $p$-torsion element of $H_{\mathfrak a}^3(R)$. Since $p\lambda_q \in
(\Delta_1^q, \ \Delta_2^q, \ \Delta_3^q)R$, it is immediate that
$$
p\cdot\eta_q = [p\lambda_q + (\Delta_1^q, \ \Delta_2^q, \ \Delta_3^q)R] = 0 \ 
\in \ H_{\mathfrak a}^3(R).
$$ 
We shall show that Hochster's conjecture is false by showing that for all 
$q=p^e$, $\eta_q=0$ in $H_{\mathfrak a}^3(R)$. More specifically we show that 
if $k=q-1$ then
$$
\lambda_q (\Delta_1 \Delta_2 \Delta_3)^{k} \in
(\Delta_1^{q+k}, \ \Delta_2^{q+k}, \ \Delta_3^{q+k})R.
$$

We record the following result, \cite[Proposition 4.1]{PS}, from which it
follows that $H_{\mathfrak a}^3(R/pR)=0$. We use this result with $A=R/P$ and
$I = (\overline{\Delta}_1, \overline{\Delta}_2, \overline{\Delta}_3)$. Then
$\height I = 2$, and it is well known that the determinantal ring $A/I$ is
Cohen-Macaulay.

\begin{prop}[Peskine-Szpiro]
Let $A$ be a regular domain of prime characteristic $p > 0$, and let $I$ be an 
ideal of $A$ such that $A/I$ is Cohen-Macaulay. If $\height I = h$, then 
$H_I^i(A)=0$ for $i > h$.
\label{PS}
\end{prop}

\begin{remark}
For a field $K$, let $R_K = K\otimes_{\mathbb Z}R$ where $R={\mathbb
Z}[u,v,w,x,y,z]$ as above. We record an argument due to Hochster which shows
that $H_{\mathfrak a}^3(R_{\mathbb Q})$ is nonzero.  Consider the action of
$SL_2({\mathbb Q})$ on $R_{\mathbb Q}$ where where $\alpha \in SL_2({\mathbb
Q})$ sends the entries of the matrix $M$ to the entries of $\alpha M$. The
invariant subring for this action is $A={\mathbb Q}[\Delta_1, \Delta_2,
\Delta_3]$ and, since $SL_2({\mathbb Q})$ is linearly reductive, $A$ is a
direct summand of $R_{\mathbb Q}$ as an $A$--module. Consequently $H_{\mathfrak
a}^3(A)$ is a nonzero submodule of $H_{\mathfrak a}^3(R_{\mathbb Q})$. Hence
while  $H_{\mathfrak a}^3(R_{{\mathbb Z}/p{\mathbb Z}})=0$, we have
$H_{\mathfrak a}^3(R_{\mathbb Q}) \neq 0$; it is then only natural to expect
that the study of $H_{\mathfrak a}^3(R)$ would be interesting!   
\end{remark}

\section{The multi-grading}

We now work with a fixed prime integer $p$, an arbitrary prime power $q=p^e$, 
and set $\lambda=\lambda_q$. We first use a multi-grading to reduce the
question whether 
$$
\lambda (\Delta_1 \Delta_2 \Delta_3)^{k} \in
(\Delta_1^{q+k}, \ \Delta_2^{q+k}, \ \Delta_3^{q+k})R \qquad (\#)
$$
to a question in a polynomial ring in three variables. We assign weights as 
follows:
\begin{align*}
u: (1,0,0,0) \qquad \qquad\qquad  & x: (1,0,0,1) \\
v: (0,1,0,0) \qquad \qquad\qquad  & y: (0,1,0,1) \\
w: (0,0,1,0) \qquad \qquad\qquad  & z: (0,0,1,1) 
\end{align*}
With this grading, 
$$
\deg(\Delta_1) = (0,1,1,1), \quad \quad
\deg(\Delta_2) = (1,0,1,1), \quad \quad
\deg(\Delta_3) = (1,1,0,1),
$$
and $\lambda$ is a homogeneous element of degree $(q,q,q,q)$. Hence in a
homogeneous equation of the form
$$
\lambda (\Delta_1\Delta_2\Delta_3)^{k}
=c_1\Delta_1^{q+k}+c_2\Delta_2^{q+k}+c_3\Delta_3^{q+k},
$$
we must have 
$$
\deg(c_1) = (q+2k,k,k,2k), \ \deg(c_2) = (k,q+2k,k,2k), \ 
\deg(c_3) = (k,k,q+2k,2k). 
$$
We use this to examine the monomials which can occur in $c_i$. If $c_1$ 
involves a monomial of the form $u^a v^b w^c x^d y^i z^j$, then
$$
(a+d, \ b+i, \ c+j, \ d+i+j) = (q+2k, \ k, \ k, \ 2k)
$$
and so $b=k-i, \ c=k-j, \ d=2k-i-j$ and $a=q+i+j$. Hence $c_1$ is a ${\mathbb
Z}$-linear combination  of monomials of the form
$$
u^{q+i+j}v^{k-i}w^{k-j}x^{2k-i-j}y^iz^j = u^q(uy)^i(vx)^{k-i}(uz)^j(wx)^{k-j}. 
$$
Let $[s_1,s_2]^k$ denote the set of all monomials $s_1^is_2^{k-i}$ for $0 \le i
\le k$ and $[s_1,s_2]\cdot[t_1,t_2]$ denote the set of products of all pairs of
monomials from $[s_1,s_2]$ and $[t_1,t_2]$. With this notation, $c_1$ is a
${\mathbb Z}$-linear combination of monomials in 
$u^q\cdot [uy,vx]^k \cdot [uz,wx]^k$. After similar computations for $c_2$ and 
$c_3$, we may conclude that $(\#)$ holds if and only if $\lambda (\Delta_1 
\Delta_2 \Delta_3)^{k}$ is a ${\mathbb Z}$-linear combination of elements from
\begin{align*}
& \Delta_1^{q+k}u^q  \cdot [uy,vx]^k  \cdot [uz,wx]^k, \\  
& \Delta_2^{q+k}v^q  \cdot [vz,wy]^k  \cdot [vx,uy]^k, \\
& \Delta_3^{q+k}w^q  \cdot [wx,uz]^k  \cdot [wy,vz]^k. 
\end{align*} 
We may divide throughout by the element $(uvw)^{q+2k}$, and study this issue in 
the polynomial ring $\mathbb{Z}[\frac{x}{u},\frac{y}{v},\frac{z}{w}]$. Let 
$$
A=\frac{z}{w}-\frac{x}{u},\qquad B=\frac{x}{u}-\frac{y}{v},\qquad T=-\frac{x}{u}.
$$
The condition $(\#)$ is then equivalent to the statement that the element
$$
\frac{\lambda (\Delta_1 \Delta_2 \Delta_3)^{k}}{(uvw)^{q+2k}} 
= \frac{1}{p} \Big((A+B)^q+(-A)^q+(-B)^q \Big) \Big( (A+B)AB \Big)^k
$$
is a ${\mathbb Z}$-linear combination, in the polynomial ring 
$\mathbb{Z}[A,B,T]$, of elements of
\begin{align*}   
&(A+B)^{q+k} \cdot [T,A]^k \cdot [T,B]^k, \\ 
& A^{q+k} \cdot [T,B]^k \cdot [T+B, A+B]^k, \\ 
& B^{q+k} \cdot [T,A]^k \cdot [T-A, A+B]^k.
\end{align*} 
We show that the above statement is indeed true for $k=q-1$ using the following 
equational identity, which will be proved in the next section:
\begin{align*}
& (A+B)^{2k+1} \sum_{n=0}^k \binom{k}{n} T^n \sum_{i=0}^n (-1)^{k+i} 
A^{k-i} B^{k-n+i} \binom{k+i}{k} \binom{k+n-i}{k} \\
& -A^{2k+1} \sum_{n=0}^k \binom{k}{n} (-1)^n (T+B)^n \sum_{i=0}^n 
(A+B)^{k-i} B^{k-n+i} \binom{k+i}{k} \binom{k+n-i}{k} \\
& -B^{2k+1} \sum_{n=0}^k \binom{k}{n} (T-A)^n \sum_{i=0}^n 
(A+B)^{k-i} A^{k-n+i} \binom{k+i}{k} \binom{k+n-i}{k} =0.
\end{align*}
Separating the terms with $n=0$ and dividing by $p$, the above identity gives us
\begin{align*}
\frac{1}{p} & \Big((A+B)^{k+1}(-1)^k-A^{k+1}-B^{k+1}\Big) \Big(AB(A+B)\Big)^k \\
 = &
-(A+B)^{2k+1} \sum_{n=1}^k \binom{k}{n} T^n \sum_{i=0}^n (-1)^{k+i} 
A^{k-i} B^{k-n+i} \frac{1}{p} \binom{k+i}{k} \binom{k+n-i}{k} \\
 & +A^{2k+1} \sum_{n=1}^k \binom{k}{n} (-1)^n (T+B)^n \sum_{i=0}^n 
(A+B)^{k-i} B^{k-n+i} \frac{1}{p}\binom{k+i}{k} \binom{k+n-i}{k} \\
 & +B^{2k+1} \sum_{n=1}^k \binom{k}{n} (T-A)^n \sum_{i=0}^n 
(A+B)^{k-i} A^{k-n+i} \frac{1}{p} \binom{k+i}{k} \binom{k+n-i}{k}.
\end{align*}
To ensure that the expression on the right hand side of the equality is indeed
a ${\mathbb Z}$-linear combination of elements of the required form, it 
suffices to establish that the prime integer $p$ divides
$$
\binom{k+i}{k} \binom{k+n-i}{k}
$$
when $1 \le n \le k$ and $0 \le i \le n$, where $k=q-1=p^e-1$.
Note that this condition ensures that $i \ge 1$ or $n-i \ge 1$, and so it
suffices to show that $p$ divides 
$$
\binom{k+r}{k} = \binom{p^e-1+r}{p^e-1}
$$ 
for $1 \le r \le p^e-1$. Since
\begin{align*}
\binom{p^e-1+r}{p^e-1} &=\left(\frac{p^e}{r}\right)\left(\frac{p^e+1}{1}\right) 
\left(\frac{p^e+2}{2}\right) \cdots \left(\frac{p^e+r-1}{r-1}\right) \\
&=\left(\frac{p^e}{r}\right)\left(\frac{p^e}{1}+1\right) 
\left(\frac{p^e}{2}+1\right) \cdots \left(\frac{p^e}{r-1}+1\right),
\end{align*}
this is easily seen to be true.

\section{The equational identity}

We first record some identities with binomial coefficients that we shall be
using. These identities can be easily proved using Zeilberger's algorithm (see
\cite{PWS}) and the Maple package {\tt EKHAD}, but we include these proofs for
the sake of completeness.

When the range of a summation is not specified, it is assumed to extend over
all integers. We set $\binom{k}{i}=0$ if $ i<0$ or if $k<i$.

\begin{lemma}
\label{iden}
\begin{align}
& \sum_{i=0}^m (-1)^i \binom{2k+1}{s+i}\binom{k+i}{k} \binom{k+m-i}{k} =
\begin{cases}
\binom{m+s-k-1}{m} \binom{2k+m+1}{m+s} \, \text{if } s > k, \\
0 \qquad  \ \text{ if } k \ge s \ge k+1-m,
\end{cases} \\
& \sum_{j=0}^s (-1)^j \binom{s}{j} \binom{a+j}{k} 
= (-1)^s \binom{a}{k-s}, \\
& \sum_{i=0}^s \binom{k-i}{k-s} \binom{k+i}{k} \binom{k+m-i}{m}
= \binom{k+m-s}{m} \binom{2k+m+1}{s}. 
\end{align}
\end{lemma}

\begin{proof}
(1) Let 
\begin{align*}
&F(m,i) = (-1)^i \binom{2k+1}{s+i}\binom{k+i}{k} \binom{k+m-i}{k}, \\
&G(m,i) = \frac{i(s+i)(k+m+1-i)}{m+1-i} F(m,i), \quad \text{ and } \quad
H(m)=\sum_i F(m,i).
\end{align*}
It is easily seen that
\begin{align*}
&G(m,i+1) = -(k+i+1)(2k+1-s-i)F(m,i) \quad \text{ and } \quad \\
&F(m+1,i) = \frac{k+m+1-i}{m+1-i} F(m,i),
\end{align*}
and these can then be used to verify the relation
$$
G(m,i+1) - G(m,i) = (2k+m+2)(m+s-k)F(m,i) - (m+1)(m+s+1)F(m+1,i). 
$$
Summing with respect to $i$ gives us
$$
0 = (2k+m+2)(m+s-k)H(m) - (m+1)(m+s+1)H(m+1). 
$$
Using this recurrence,
$$
H(m) = \frac{(2k+m+1)\cdots (2k+2)(m+s-k-1)\cdots (s-k)}
{(m)\cdots (1) (m+s)\cdots (s+1)}H(0),
$$
and the required result follows.

(2) Let 
$$
F(s,j) = (-1)^j \binom{s}{j} \binom{a+j}{k}, \
G(s,j) = \frac{-j(a+j-k)}{s+1-j} F(s,j),
$$
and $H(s)=\sum_j F(s,j)$. It is easily seen that
\begin{align*}
&G(s,j+1) = (a+j+1)F(s,j) \quad \text{ and } \quad 
F(s+1,j) = \frac{s+1}{s+1-j} F(s,j),
\end{align*}
and these can then be used to verify the relation
$$
G(s,j+1) - G(s,j) = (k-s)F(s,j) + (a-k+s+1)F(s+1,j). 
$$
Summing with respect to $j$ gives us
$$
0 = (k-s)H(s) + (a-k+s+1)H(s+1).
$$
Using this recurrence,
\begin{align*}
H(s) = \frac{(-1)^s (k-s+1)\cdots (k)}{(a-k+s)\cdots (a-k+1)} H(0)
= (-1)^s \binom{a}{k-s}.
\end{align*}

(3) Let 
$$
F(s,i) = \binom{k-i}{k-s} \binom{k+i}{k} \binom{k+m-i}{m}, \
G(s,i) = \frac{i(k-s)(k+m+1-i)}{s+1-i} F(s,i),
$$
and $H(s)=\sum_i F(s,i)$. It is easily seen that
$$
G(s,i+1) = (k-s) (k+i+1) F(s,i) \quad \text{ and } \quad 
F(s+1,i) = \frac{k-s}{s+1-i} F(s,i),
$$
and these can then be used to verify the relation
$$
G(s,i+1) - G(s,i) = (k-s)(2k+m-s+1)F(s,i) - (s+1)(k+m-s)F(s+1,i). 
$$
Summing with respect to $j$ gives us
$$
0 = (k-s)(2k+m-s+1)H(s) - (s+1)(k+m-s)H(s+1). 
$$
Using this recurrence,
\begin{align*}
H(s) & = \frac{(k-s+1)\cdots (k)(2k+m-s+2)\cdots (2k+m+1)}
{(s)\cdots (1) (k+m-s+1)\cdots (k+m)}H(0) \\
& = \binom{k+m-s}{m} \binom{2k+m+1}{s}.
\end{align*}
\end{proof}

We are now ready to prove the equational identity
\begin{align*}
&(A+B)^{2k+1} \sum_{n=0}^k \binom{k}{n} T^n \sum_{i=0}^n (-1)^{k+i} 
A^{k-i} B^{k-n+i} \binom{k+i}{k} \binom{k+n-i}{k} \\
&-A^{2k+1} \sum_{n=0}^k \binom{k}{n} (-1)^n (T+B)^n \sum_{i=0}^n 
(A+B)^{k-i} B^{k-n+i} \binom{k+i}{k} \binom{k+n-i}{k} \\
&-B^{2k+1} \sum_{n=0}^k \binom{k}{n} (T-A)^n \sum_{i=0}^n 
(A+B)^{k-i} A^{k-n+i} \binom{k+i}{k} \binom{k+n-i}{k} =0.
\end{align*}
Examining the coefficient of $T^m$ for all $0 \le m \le k$, we need to show
\begin{align*}
& (A+B)^{2k+1}\binom{k}{m} \sum_{i=0}^m (-1)^{k+i} 
A^{k-i} B^{k-m+i} \binom{k+i}{k} \binom{k+m-i}{k} \\
& -A^{2k+1} \sum_{n=m}^k \binom{k}{n} \binom{n}{m} (-1)^n \sum_{i=0}^n 
(A+B)^{k-i} B^{k-m+i} \binom{k+i}{k} \binom{k+n-i}{k} \\
& -B^{2k+1} \sum_{n=m}^k \binom{k}{n} \binom{n}{m} (-1)^{n-m} \sum_{i=0}^n 
(A+B)^{k-i} A^{k-m+i} \binom{k+i}{k} \binom{k+n-i}{k} =0.
\end{align*}
Dividing by $\binom{k}{m}(AB)^{k-m}$ and using the fact that
$$
\binom{k}{m}\binom{k-m}{n-m} = \binom{k}{n}\binom{n}{m},
$$
we need to establish that
\begin{align*}
& (A+B)^{2k+1} \sum_{i=0}^m (-1)^{k+i} 
A^{m-i} B^i \binom{k+i}{k} \binom{k+m-i}{k} \\
& -A^{k+m+1} \sum_{n=m}^k \binom{k-m}{n-m} (-1)^n \sum_{i=0}^n 
(A+B)^{k-i} B^i \binom{k+i}{k} \binom{k+n-i}{k} \\
& -B^{k+m+1} \sum_{n=m}^k \binom{k-m}{n-m} (-1)^{n-m} \sum_{i=0}^n 
(A+B)^{k-i} A^i \binom{k+i}{k} \binom{k+n-i}{k} =0.
\end{align*}
For $0 \le r \le m-1$, the coefficient of $A^{k+m-r}B^{k+1+r}$ is
$$
\sum_{i=0}^{m}(-1)^{k+i} \binom{2k+1}{k-r+i} \binom{k+i}{k} \binom{k+m-i}{k}
$$
which is zero by lemma \ref{iden} (1) since $k \ge k-r \ge k+1-m$. 

For $0 \le r \le k$, the coefficient  of $A^{k+m+1+r}B^{k-r}$ as well as
the coefficient of  $(-1)^m A^{k-r} B^{k+m+1+r}$ is 
\begin{align*}
&\sum_{i=0}^m (-1)^{k+i} \binom{2k+1}{k+1+r+i}\binom{k+i}{k} \binom{k+m-i}{k} \\
&\qquad \qquad\qquad  - \sum_{n=m}^k (-1)^n \binom{k-m}{n-m} 
\sum_i \binom{k-i}{r} \binom{k+i}{k} \binom{k+n-i}{k} \\
&=(-1)^k\binom{m+r}{m} \binom{2k+m+1}{k-r} \\
&\qquad \qquad\qquad  -\sum_{i=0}^{k-r} \binom{k-i}{r} \binom{k+i}{k} 
\sum_{n=m}^k (-1)^n \binom{k-m}{n-m} \binom{k+n-i}{k} \\
&=(-1)^k\binom{m+r}{m} \binom{2k+m+1}{k-r} 
- \sum_{i=0}^{k-r} (-1)^k \binom{k-i}{r} \binom{k+i}{k} \binom{k+m-i}{m} \\
&=0,
\end{align*}
using identities established in lemma 3.1.

\section{A local cohomology module with infinitely many associated prime
ideals}\label{elem}

Consider the hypersurface of mixed characteristic
$$
R={\mathbb Z}[U,V,W,X,Y,Z]/(UX+VY+WZ)
$$
and the ideal ${\mathfrak a}=(x,y,z)R$. We show that the local cohomology
module $H_{\mathfrak a}^3(R)$ has $p$-torsion elements for infinitely many
prime integers $p$, and consequently that it has infinitely many associated
prime ideals. We have
$$
H_{\mathfrak a}^3(R) = \varinjlim \frac{R}{(x^k, y^k, z^k)R}
$$ 
where the maps in the direct limit system are induced by multiplication by the 
element $xyz$. Let $p$ be a prime integer. It is easily seen that 
$$
\lambda=\frac{(ux)^p + (vy)^p + (wz)^p}{p}
$$
has integer coefficients, and we claim that the element
$$
\eta = [\lambda + (x^p, \ y^p, \ z^p)R] \ \in \ H_{\mathfrak a}^3(R)
$$
is nonzero and $p$-torsion. It is immediate that $p\cdot\eta = 0$, and what 
really needs to be established is that $\eta$ is a nonzero element of 
$H_{\mathfrak a}^3(R)$, i.e., that
$$
\lambda (xyz)^k \notin 
(x^{p+k}, \ y^{p+k}, \ z^{p+k})R \ \text { for all } \ k \in {\mathbb N}. 
$$
We shall accomplish this using an ${\mathbb N}^4$-grading. We assign weights as 
follows:
\begin{align*}
u: (0,1,1,1)  \qquad \qquad\qquad  & x: (1,0,0,0) \\
v: (1,0,1,1)  \qquad \qquad\qquad  & y: (0,1,0,0) \\
w: (1,1,0,1)  \qquad \qquad\qquad  & z: (0,0,1,0) 
\end{align*}
With this grading, $\lambda$ is a homogeneous element of degree $(p,p,p,p)$. 
Hence in a homogeneous equation of the form
$$
\lambda (xyz)^{k}
=c_1x^{p+k}+c_2y^{p+k}+c_3z^{p+k},
$$
we must have 
$$
\deg(c_1) = (0,p+k,p+k,p), \ \deg(c_2) = (p+k,0,p+k,p), \ 
\deg(c_3) = (p+k,p+k,0,p). 
$$
We may use this to examine the monomials which can occur in $c_i$, and it is 
easily seen that the only monomial that can occur in $c_1$ with a nonzero
coefficient is $u^p y^k z^k$, and 
similarly that $c_2$ is an integer multiple of $v^p z^k x^k$ and $c_3$ 
is an integer multiple of $w^p x^k y^k$. Hence 
$\lambda (xyz)^k \in (x^{p+k}, y^{p+k}, z^{p+k})R$ if and only if
\begin{align*}
\lambda (xyz)^k 
&\in (u^p y^k z^k x^{p+k}, \ v^p z^k x^k y^{p+k}, \ w^p x^k y^k z^{p+k})R \\
&= (xyz)^k\Big((ux)^p, \ (vy)^p, \ (wz)^p\Big)R,
\end{align*}
i.e., if and only if $\lambda \in ((ux)^p, (vy)^p, (wz)^p)R$. To complete our 
argument it suffices to show that $\lambda \notin (p,(ux)^p, (vy)^p, (wz)^p)R$, 
i.e., that  
$$
\frac{(ux)^p+(vy)^p+(-1)^p(ux+vy)^p}{p} \notin \Big(p,\ (ux)^p,\ (vy)^p \Big)R.
$$
After making the specializations $u \mapsto 1, \ v \mapsto 1, \ w \mapsto 1$,
it is enough to verify that
$$
\frac{x^p+y^p+(-1)^p(x+y)^p}{p} \notin (p, \ x^p, \ y^p){\mathbb Z}[x,y]
$$
and this holds since the coefficient of $x^{p-1}y$ in 
$(x^p+y^p+(-1)^p(x+y)^p)/p$ is $(-1)^p$.

\section*{Acknowledgments}

The use of the Maple package {\tt EKHAD} is gratefully acknowledged. I would 
like to thank Mel Hochster for first bringing the conjecture to my attention,
and I am also grateful to Ian Aberbach and Reinhold H\"ubl for pointing out an
error in an earlier version of this manuscript.


\begin{thebibliography}{PWS}

\bibitem[BL]{BL} M. P. Brodmann and A. Lashgari Faghani, {\em A finiteness
result for associated primes of local cohomology modules}, Proc. Amer. Math.
Soc., to appear.

\bibitem[BRS]{BRS} M. P. Brodmann, C. Rotthaus and R. Y. Sharp, {\em On
annihilators and associated primes of local cohomology modules}, J. Pure Appl.
Alg., to appear.

\bibitem[He]{He} M. Hellus, {\em On the set of associated primes of a local
cohomology module}, preprint.

\bibitem[Hu]{Hu} C. Huneke, {\em Problems on local cohomology}, in: Free
resolutions in commutative algebra and algebraic geometry (Sundance, Utah,
1990), pp. 93--108, Research Notes in Mathematics, 2, Jones and Bartlett
Publishers, Boston, MA, 1992. 

\bibitem[HS]{HS} C. Huneke and R. Sharp, {\em Bass numbers of local cohomology
modules}, Trans. Amer. Math. Soc. {\bf 339} (1993), 765--779. 

\bibitem[Ly1]{Ly1} G. Lyubeznik, {\em Finiteness properties of local cohomology
modules (an application of $D$-modules to commutative algebra)}, Invent. Math.
{\bf 113} (1993), 41--55.

\bibitem[Ly2]{Ly2} G. Lyubeznik, {\em $F$-modules: applications to local
cohomology and $D$-modules in characteristic $p>0$}, J. Reine Angew. Math. {\bf
491} (1997), 65--130.

\bibitem[Ly3]{Ly3} G. Lyubeznik, {\em Finiteness properties of local cohomology
modules: a characteristic-free approach}, J. Pure Appl. Alg., to appear.

\bibitem[Ly4]{Ly4} G. Lyubeznik, {\em Finiteness properties of local cohomology
modules for regular local rings of mixed characteristic: the unramified case},
Comm. Alg., to appear.

\bibitem[PS]{PS} C. Peskine and L. Szpiro, {\em Dimension projective finie et
cohomologie locale}, Inst. Hautes \'Etudes Sci. Publ. Math. {\bf 42} (1973),
47--119.

\bibitem[PWZ]{PWS} M. Petkov\v sek, H. S. Wilf and D. Zeilberger, $A=B$, With a
foreword by Donald E. Knuth, A K Peters, Ltd., Wellesley, MA, 1996. 

\end{thebibliography}
\end{document}